\input amstex
\input amsppt.sty
\magnification=\magstep1
\hsize=32truecc
\vsize=22.2truecm
\baselineskip=16truept
\NoBlackBoxes
\TagsOnRight \pageno=1 \nologo
\def\Z{\Bbb Z}
\def\N{\Bbb N}

\def\Q{\Bbb Q}

\def\C{\Bbb C}
\def\l{\left}
\def\r{\right}
\def\bg{\bigg}
\def\({\bg(}
\def\[{\bg\lfloor}
\def\){\bg)}
\def\]{\bg\rfloor}
\def\t{\text}
\def\f{\frac}

\def\se {\subseteq}

\def\bi{\binom}
\def\eq{\equiv}

\def\ls{\leqslant}
\def\gs{\geqslant}

\def\Tor{\roman{Tor}}
\def\sign{\roman{sign}}

\def\per{\roman{per}}

\def\Proof{\noindent{\it Proof}}

\def\Remark{\medskip\noindent{\it  Remark}}

\hbox{J. Algebraic Combin., in press.}
\bigskip
\topmatter
\title On permutations of $\{1,\ldots,n\}$ and related topics\endtitle
\author Zhi-Wei Sun\endauthor
\leftheadtext{Zhi-Wei Sun}
\rightheadtext{On permutations of $\{1,\ldots,n\}$ and related topics}
\affil Department of Mathematics, Nanjing University\\
 Nanjing 210093, People's Republic of China
  \\  zwsun\@nju.edu.cn
  \\ {\tt http://maths.nju.edu.cn/$\sim$zwsun}
\endaffil
\abstract In this paper we study combinatorial aspects of permutations of $\{1,\ldots,n\}$
and related topics. In particular, we prove that there is a unique permutation $\pi$ of $\{1,\ldots,n\}$
such that all the numbers $k+\pi(k)$ ($k=1,\ldots,n$) are powers of two.
We also show that $n\mid\per[i^{j-1}]_{1\ls i,j\ls n}$ for any integer $n>2$.
We conjecture that if a group $G$ contains no element of order among $2,\ldots,n+1$ then
any $A\se G$ with $|A|=n$ can be written as  $\{a_1,\ldots,a_n\}$
with $a_1,a_2^2,\ldots,a_n^n$ pairwise distinct. This conjecture is confirmed when $G$ is a torsion-free abelian group. We also prove that for any finite subset $A$ of a torsion-free abelian group $G$ with $|A|=n>3$, there is a numbering $a_1,\ldots,a_n$ of all the elements of $A$ such that all the $n$ sums
$$a_1+a_2+a_3,\ a_2+a_3+a_4,\ \ldots,\ a_{n-2}+a_{n-1}+a_n,\ a_{n-1}+a_n+a_1,\ a_n+a_1+a_2$$
are pairwise distinct, and conjecture that this remains valid if $G$ is cyclic.
\endabstract
\thanks 2020 {\it Mathematics Subject Classification}. Primary 05A05, 11B75;
Secondary 11B13, 11B39, 20D60.
\newline  \indent {\it Keywords}. Additive combinatorics, permutations, powers of two, permanents, groups.
\newline \indent Supported by the National Natural Science
Foundation of China (grant no. 11971222).
\endthanks
\endtopmatter
\document

\heading{1. Introduction}\endheading

As usual, for $n\in\Z^+=\{1,2,3,\ldots\}$ we let $S_n$ denote the symmetric group of all the permutation of $\{1,\ldots,n\}$.

Let $A=[a_{ij}]_{1\ls i,j\ls n}$ be a $(0,1)$-matrix (i.e., $a_{ij}\in\{0,1\}$ for all $i,j=1,\ldots,n$).
Then the permanent of $A$ given by
$$\per(A)=\sum_{\pi\in S_n}a_{1\pi(1)}\cdots a_{n\pi(n)}$$
is just the number of permutations $\pi\in S_n$ with $a_{k\pi(k)}=1$ for all $k=1,\ldots,n$.

In 2002, B. Cloitre proposed the sequence [Cl, A073364] on OEIS whose $n$-th term $a(n)$ is the number of permutations
$\pi\in S_n$ with $k+\pi(k)$ prime for all $k=1,\ldots,n$. Clearly, $a(n)=\per(A)$, where $A$ is a matrix of order $n$ whose $(i,j)$-entry ($1\ls i,j\ls n$) is $1$ or $0$ according as $i+j$ is prime or not.
In 2018 P. Bradley [Br] proved that $a(n)>0$ for all $n\in\Z^+$.

Our first theorem is an extension of Bradley's result.

\proclaim{Theorem 1.1} Let $(a_1,a_2,\ldots)$ be an integer sequence with $a_1=2$ and $a_k<a_{k+1}\ls 2a_{k}$
for all $k=1,2,3\ldots$. Then, for any positive integer $n$, there exists a permutation $\pi\in S_n$
with $\pi^2=I_n$ such that
$$\{k+\pi(k):\ k=1,\ldots,n\}\se\{a_1,a_2,\ldots\},\tag1.1$$
where $I_n$ is the identity of $S_n$ with $I_n(k)=k$ for all $k=1,\ldots,n$.
\endproclaim

Recall that the Fiboncci numbers $F_0,F_1,\ldots$ and the Lucas numbers $L_0,L_1,\ldots$ are defined by
$$F_0=0,\ F_1=1,\ \t{and}\ F_{n+1}=F_n+F_{n-1}\ (n=1,2,3,\ldots),$$
and
$$L_0=2,\ L_1=1,\ \t{and}\ L_{n+1}=L_n+L_{n-1}\ (n=1,2,3,\ldots).$$
If we apply Theorem 1.1 with the sequence $(a_1,a_2,\ldots)$ equal to $(F_3,F_4,\ldots)$
or $(L_0,L_2,L_3,\ldots)$, then we immediately obtain the following consequence.

\proclaim{Corollary 1.1} Let $n\in\Z^+$. Then there is a permutation $\sigma\in S_n$
with $\sigma^2=I_n$ such that
all the sums $k+\sigma(k)\ (k=1,\ldots,n)$ are Fibonacci numbers. Also,
there is a permutation $\tau\in S_n$ with $\tau^2=I_n$ such that
all the numbers $k+\tau(k)\ (k=1,\ldots,n)$ are Lucas numbers.
\endproclaim
\Remark\ 1.1. Let $f(n)$ be the number of permutations $\sigma\in S_n$ such that all the sums $k+\sigma(k)\ (k=1,\ldots,n)$ are Fibonacci numbers. Via {\tt Mathematica} we find that
$$(f(1),\ldots,f(20))=(1,1,1,2,1,2,4,2,1,4,4,20,4,5,1,20,24,8,96,200).$$
For example, $\pi=(2,3)(4,9)(5,8)(6,7)$ is the unique permutation in $S_9$ such that
all the numbers $k+\pi(k)\ (k=1,\ldots,9)$ are Fibonacci numbers.
\medskip

Recall that those integers $T_n=n(n+1)/2\ (n=0,1,2,\ldots)$ are called triangular numbers.
Note that $T_n-T_{n-1}=n\ls T_{n-1}$ for every $n=3,4,\ldots$.
Applying Theorem 1.1 with $(a_1,a_2,a_3,\ldots)=(2,T_2,T_3,\ldots)$, we immediately get the following corollary.

\proclaim{Corollary 1.2} For any $n\in\Z^+$, there is a permutation $\pi\in S_n$
with $\pi^2=I_n$ such that
each of the sums $k+\pi(k)\ (k=1,\ldots,n)$ is either $2$ or a triangular number.
\endproclaim
\Remark\ 1.2. When $n=4$, we may take $\pi=(2,4)$ to meet the requirement in Corollary 1.2.
Note that $1+1=3=T_2$ and $2+4=3+3=T_3$.
\medskip

Our next theorem focuses on permutations involving powers of two.

\proclaim{Theorem 1.2} Let $n$ be any positive integer. Then there is a unique permutation $\pi_n\in S_n$
such that all the numbers $k+\pi_n(k)\ (k=1,\ldots,n)$ are powers of two. In other words, for the $n\times n$
matrix $A$ whose $(i,j)$-entry is $1$ or $0$ according as $i+j$ is a power of two or not, we have $\per(A)=1$.
\endproclaim
\Remark\ 1.3. Note that the number of 1's in the matrix $A$ given in Theorem 1.2 coincides with
$$\sum_{k=0}^{\lfloor \log_2n\rfloor+1}\sum_{1\ls i,j\ls n\atop\ i+j=2^k}1=\sum_{k=0}^{\lfloor\log_2n\rfloor}(2^k-1)
+\sum_{i=2^{\lfloor\log_2n\rfloor+1}-n}^n1=2n-\lfloor\log_2n\rfloor-1.$$

{\it Example} 1.1. Here we list $\pi_n$ in Theorem 1.2 for $n=1,\ldots,11$:
$$\gather \pi_1=(1),\ \pi_2=(1),\ \pi_3=(1,3),\ \pi_4=(1,3),\ \pi_5=(3,5),\ \pi_6=(2,6)(3,5),
\\\pi_7=(1,7)(2,6)(3,5),\ \pi_8=(1,7)(2,6)(3,5),\ \pi_9=(2,6)(3,5)(7,9),
\\\pi_{10}=(3,5)(6,10)(7,9),\ \pi_{11}=(1,3)(5,11)(6,10)(7,9).
\endgather$$

Theorem 1.2 has the following consequence.

\proclaim{Corollary 1.3} For any $n\in\Z^+$, there is a unique permutation $\pi\in S_{2n}$
such that $k+\pi(k)\in\{2^a-1:\ a\in\Z^+\}$ for all $k=1,\ldots,2n$.
\endproclaim

Now we turn to our results of new types.

\proclaim{Theorem 1.3}
{\rm (i)} Let $p$ be any odd prime. Then there is no $\pi\in S_{p-1}$ such that
all the $p-1$ numbers $k\pi(k)\ (k=1,\ldots,p-1)$ are pairwise incongruent modulo $p$.
Also,
$$\per[i^{j-1}]_{1\ls i,j\ls p-1}\eq0\pmod p.\tag1.2$$

{\rm (ii)} We have
$$\per[i^{j-1}]_{1\ls i,j\ls n}\eq0\pmod n\ \ \t{for all}\ n=3,4,5,\ldots.\tag1.3$$
\endproclaim
\Remark\ 1.4. In contrast with Theorem 1.3, it is well-known that
$$\det[i^{j-1}]_{1\ls i,j\ls n}=\prod_{1\ls i<j\ls n}(j-i)=1!2!\ldots (n-1)!$$
and in particular
$$\det[i^{j-1}]_{1\ls i,j\ls p-1},\,\det[i^{j-1}]_{1\ls i,j\ls p}\not\eq0\pmod p$$
for any odd prime $p$.
\medskip

In additive combinatorics, there are some interesting topics involving both permutations and finite abelian groups, see. e.g., [FSX] and [GS].
Below we present two novel theorems on permutations involving groups.

\proclaim{Theorem 1.4} {\rm (i)} Let $a_1,\ldots,a_n$ be distinct elements of a torsion-free abelian group $G$. Then there is a permutation
$\pi\in S_n$ such that all those $ka_{\pi(k)}\ (k=1,\ldots,n)$ are pairwise distinct.

{\rm (ii)} Let $a,b,c$ be three distinct elements of a group $G$ such that none of them has order $2$ or $3$.
Then $a^{\sigma(1)}$ and $b^{\sigma(2)}$ are distinct for some $\sigma\in S_2$. Also,
$a^{\tau(1)},b^{\tau(2)},c^{\tau(3)}$ are pairwise distinct for some $\tau\in S_3$.
\endproclaim
\Remark\ 1.5. On the basis of this theorem, we will formulate a general conjecture for groups in Section 4.
\medskip

\proclaim{Theorem 1.5}  For any $n>3$ distinct elements $a_1,a_2,\ldots,a_n$ of a torsion-free abelian group $G$, there is a permutation $\pi\in S_n$ such that all the $n$ sums
$$b_1+b_2+b_3,\ b_2+b_3+b_4,\ \ldots,\ b_{n-2}+b_{n-1}+b_n,\ b_{n-1}+b_n+b_1,\ b_n+b_1+b_2$$
are pairwise distinct, where $b_k=a_{\pi(k)}$ for $k=1,\ldots,n$.
\endproclaim
\Remark\ 1.6.  By Remark 1.2 of Sun [S19], for any finite subset $A$ of a torsion-free abelian group with $|A|=n>2$
we may write $A$ as $\{a_1,\ldots,a_n\}$ such that $a_1+a_2,\ldots,a_{n-1}+a_n,a_n+a_1$ are pairwise distinct.
\medskip

We are going to prove Theorems 1.1-1.3 and Corollary 1.3 in the next section, and
show Theorems 1.4-1.5 in Section 3. We will pose some conjectures in Section 4.

\heading{2. Proofs of Theorems 1.1-1.3 and Corollary 1.3}\endheading

\medskip
\noindent{\it Proof of Theorem 1.1}. For convenience, we set $a_0=1$ and $A=\{a_1,a_2,a_3,\ldots\}$.
We use induction on $n\in\Z^+$ to show the desired result.

For $n=1$, we take $\pi(1)=1$ and note that $1+\pi(1)=2=a_1\in A$.

Now let $n\gs2$ and assume the desired result for smaller values of $n$.
Choose $k\in\N$ with $a_k\ls n<a_{k+1}$, and write $m=a_{k+1}-n$. Then
$1\ls m\ls 2a_k-n\ls2n-n=n$. Let $\pi(j)=a_{k+1}-j$ for $j=m,\ldots,n$.
Then  $$\{\pi(j):\ j=m,\ldots,n\}=\{m,\ldots,n\},$$
and $\pi(\pi(j))=j$ for all $j=m,\ldots,n$.

{\it Case} 1. $m=1$.

In this case, $\pi\in S_n$ and $\pi^2=I_n$.

{\it Case} 2. $m=n$.

In this case, $a_{k+1}=2n\gs 2a_k$. On the other hand, $a_{k+1}\ls 2a_k$. So, $a_{k+1}=2a_k$ and $a_k=n$.
Let $\pi(j)=n-j=a_k-j$ for all $0<j<n$. Then $\pi\in S_n$ and $j+\pi(j)\in\{a_k,a_{k+1}\}$ for all $j=1,\ldots,n$.
Note that $\pi^2(k)=k$ for all $k=1,\ldots,n$.

{\it Case} 3. $1<m<n$.

In this case, by the induction hypothesis, for some $\sigma\in S_{m-1}$
with $\sigma^2=I_{m-1}$, we have $i+\sigma(i)\in A$ for all $i=1,\ldots,m-1$. Let $\pi(i)=\sigma(i)$ for all $i=1,\ldots,m-1$. Then $\pi\in S_n$ and it meets our requirement.

In view of the above, we have completed the induction proof. \qed

\medskip
\noindent{\it Proof of Theorem 1.2}. Applying Theorem 1.1 with $a_k=2^k$ for all $k\in\Z^+$, we see that
for some $\pi\in S_n$ with $\pi^2=I_n$ all the numbers $k+\pi(k)\ (k=1,\ldots,n)$ are powers of two.

Below we use induction on $n$ to prove that the number of $\pi\in S_n$ with
$$\{k+\pi(k):\ k=1,\ldots,n\}\se\{2^a:\ a\in\Z^+\}$$
is exactly one.

The case $n=1$ is trivial.

Now let $n>1$ and assume that for each $m=1,\ldots,n-1$ there is a unique $\pi_m\in S_m$
such that all the numbers $k+\pi_m(k)\ (k=1,\ldots,m)$ are powers of two. Choose $a\in\Z^+$
with $2^{a-1}\ls n<2^a$, and write $m=2^a-n$. Then $1\ls m\ls n$.

Suppose that $\pi\in S_n$ and all the numbers $k+\pi(k)\ (k=1,\ldots,n)$ are powers of two.
If $2^{a-1}\ls k\ls n$, then
$$2^{a-1}<k+\pi(k)\ls k+n\ls 2n<2^{a+1}$$
and hence $\pi(k)=2^a-k$ since $k+\pi(k)$ is a power of two. Thus
$$\{\pi(k):\ k=2^{a-1},\ldots,n\}=\{m,\ldots,2^{a-1}\}.$$
If $k\in\{1,\ldots,2^{a-1}-1\}$ and $2^{a-1}<\pi(k)\ls n$, then
$$2^{a-1}<k+\pi(k)\ls n+n<2^{a+1},$$
hence $k+\pi(k)=2^a=m+n$ and thus $m\ls k<2^{a-1}$. So we have
$$\{\pi^{-1}(j):\ 2^{a-1}<j\ls n\}=\{m,\ldots,2^{a-1}-1\}.$$
(Note that $n-2^{a-1}=2^a-m-2^{a-1}=2^{a-1}-m$.)

By the above analysis, $\pi(k)=2^a-k$ for all $k=m,\ldots,n$, and
$$\{\pi(k):\ k=m,\ldots,n\}=\{m,\ldots,n\}.$$
Thus $\pi$ is uniquely determined if $m=1$.

Now assume that $m>1$. As $\pi\in S_n$, we must have
$$\{\pi(k):\ k=1,\ldots,m-1\}=\{1,\ldots,m-1\}.$$
Since $k+\pi(k)$ is a power of two for every $k=1,\ldots,m-1$,
by the induction hypothesis we have $\pi(k)=\pi_m(k)$ for all $k=1,\ldots,m-1$.
Thus $\pi$ is indeed uniquely determined.

In view of the above, the proof of Theorem 1.2 is now complete. \qed

\medskip
\noindent{\it Proof of Corollary 1.3}. Clearly,
$\pi\in S_{2n}$ and $k+\pi(k)\in\{2^a-1:\ a\in\Z^+\}$ for all $k=1,\ldots,2n$, if and only if
there are $\sigma,\tau\in S_n$ with $\pi(2k)=2\sigma(k)-1$ and $\pi(2k-1)=2\tau(k)$
for all $k=1,\ldots,n$ such that $k+\sigma(k),k+\tau(k)\in\{2^{a-1}:\ a\in\Z^+\}$
for all $k=1,\ldots,n$. Thus we get the desired result by applying Theorem 1.2. \qed

\proclaim{Lemma 2.1 {\rm (Alon's Combinatorial Nullstellensatz [A])}} Let $A_1,\ldots,A_n$ be finite subsets
of a field $F$ with $|A_i|>k_i$ for $i=1,\ldots,n$
where $k_1,\ldots,k_n\in\{0,1,2,\ldots\}$.
 If the coefficient
of the monomial $x_1^{k_1}\cdots x_n^{k_n}$ in $P(x_1,\ldots,x_n)\in F[x_1,\ldots,x_n]$
is nonzero and $k_1+\cdots+k_n$ is
the total degree of $P$,
then there are $a_1\in A_1,\ldots,a_n\in A_n$ such that
$P(a_1,\ldots,a_n)\not=0$.
\endproclaim

\proclaim{Lemma 2.2} Let $a_1,\ldots,a_n$ be elements of a field $F$. Then
the coefficient of $x_1^{n-1}\ldots x_n^{n-1}$ in the polynomial
$$\prod_{1\ls i<j\ls n}(x_j-x_i)(a_jx_j-a_ix_i)\in F[x_1,\ldots,x_n]$$
is $(-1)^{n(n-1)/2}\per[a_i^{j-1}]_{1\ls i,j\ls n}$.
\endproclaim
\Proof. This is easy. In fact,
 $$\align&\prod_{1\ls i<j\ls
 n}(x_j-x_i)(a_jx_j-a_ix_i)
 \\=&(-1)^{\bi n2}\det[x_i^{n-j}]_{1\ls i,j\ls n}
 \times \det[a_i^{j-1}x_i^{j-1}]_{1\ls i,j\ls n}
 \\=&(-1)^{\bi n2}\sum_{\sigma\in
 S_n}\sign(\sigma)\prod_{i=1}^nx_i^{n-\sigma(i)}
 \sum_{\tau\in
 S_n}\sign(\tau)\prod_{i=1}^na_i^{\tau(i)-1}x_i^{\tau(i)-1}.
 \endalign$$
 Therefore the coefficient of $x_1^{n-1}\ldots x_n^{n-1}$ in this polynomial is
 $$(-1)^{\bi n2}\sum_{\sigma\in S_n}\sign(\sigma)^2\prod_{i=1}^na_i^{\sigma(i)-1}
 =(-1)^{n(n-1)/2}\per[a_i^{j-1}]_{1\ls i,j\ls n}.$$
 This concludes the proof. \qed

 \Remark\ 2.1. See [DKSS] and [S08, Lemma 2.2] for similar identities and arguments.

\medskip
\noindent{\it Proof of Theorem 1.3}. (i) Let $g$ be a primitive root modulo $p$. Then,
there is a permutation $\pi\in S_{p-1}$ such that the numbers $k\pi(k)\ (k=1,\ldots,p-1)$
are pairwise incongruent modulo $p$, if and only if there is a permutation $\rho\in S_{p-1}$
such that $g^{i+\rho(i)}\ (i=1,\ldots,p-1)$ are pairwise incongruent modulo $p$
(i.e., the numbers $i+\rho(i)\ (i=1,\ldots,p-1)$ are pairwise incongruent modulo $p-1$).

Suppose that $\rho\in S_{p-1}$ and all the numbers $i+\rho(i)\ (i=1,\ldots,p-1)$ are pairwise incongruent modulo $p-1$. Then
$$\sum_{i=1}^{p-1}(i+\rho(i))\eq\sum_{j=1}^{p-1}j\pmod {p-1},$$
and hence $\sum_{i=1}^{p-1}i=p(p-1)/2\eq0\pmod{p-1}$ which is impossible.
This contradiction proves the first assertion in Theorem 1.3(i).

Now we turn to prove the second assertion in Theorem 1.3(i). Suppose that
$\per[i^{j-1}]_{1\ls i,j\ls p-1}\not\eq0\pmod p$. Then, by Lemma 2.2,
 the coefficient of $x_1^{p-2}\ldots x_{p-1}^{p-2}$
in the polynomial
$$\prod_{1\ls i<j\ls p-1}(x_j-x_i)(jx_j-ix_i)$$
is not congruent to zero modulo $p$.  Applying Lemma 2.1 with $F=\Z/p\Z$
and $A=\{k+p\Z:\ k=1,\ldots,p-1\}$, we see that
there is a permutation $\pi\in S_{p-1}$ such that all those $k\pi(k)\ (k=1,\ldots,p-1)$ are pairwise
incongruent modulo $p$, which contradicts the first assertion of Theorem 1.3(i) we have just proved.

(ii) Let $n>2$ be an integer. Then
$$\align\per[i^{j-1}]_{1\ls i,j\ls n}=&\sum_{\sigma\in S_n}\prod_{k=1}^n k^{\sigma(k)-1}
\\\eq&\sum_{\sigma\in S_n\atop\sigma(n)=1}(n-1)!\prod_{k=1}^{n-1}k^{\sigma(k)-2}
=(n-1)!\sum_{\tau\in S_{n-1}}\prod_{k=1}^{n-1}k^{\tau(k)-1}
\\=&(n-1)!\,\per[i^{j-1}]_{1\ls i,j\ls n-1}
\pmod n.\endalign$$
We want to prove that $n\mid\per[i^{j-1}]_{1\le i,j\le n}$. This holds when
$n$ is an odd prime $p$, because
$p\mid \per[i^{j-1}]_{1\ls i,j\ls p-1}$ by Theorem 1.3(i).
For $n=4$, we have
 $$\align \per[i^{j-1}]_{1\ls i,j\ls 4}\eq&3!\sum_{\tau\in S_{3}}1^{\tau(1)-1}2^{\tau(2)-1}3^{\tau(3)-1}
 \\\eq&6\l(1^{2-1}2^{1-1}3^{3-1}+1^{3-1}2^{1-1}3^{2-1}\r)\eq0\pmod4.
 \endalign$$

 Now assume that $n>4$ is composite. By the above, it suffices to show that $(n-1)!\eq0\pmod n$.
 Let $p$ be the smallest prime divisor of $n$. Then $n=pq$ for some integer $q\gs p$.
 If $p<q$, then $n=pq$ divides $(n-1)!$. If $q=p$, then $p^2=n>4$ and hence $2p<p^2$,
 thus $2n=p(2p)$ divides $(n-1)!$.

 In view of the above, we have completed the proof of Theorem 1.3. \qed

 \heading{3. Proofs of Theorems 1.4 and 1.5}\endheading

 \medskip
 \noindent{\it Proof of Theorem 1.4}. (i) The case $n=1$ is trivial. Below we let $n>1$.
 Note that the subgroup $H$ of $G$ generated by $a_1,\ldots,a_n$
 is infinite, finitely generated and torsion-free.
 Thus $H$ is isomorphic to $\Z^r$ for some positive integer $r$.
 By algebraic number theory (cf. [He]), we may take an algebraic number field $K$
 with $[K:\Q]=r$ and hence $H$ is isomorphic to the additive group $O_K$ of algebraic integers in $K$.
 Thus, without any loss of generality, we may simply assume that $G$ is the additive group $\C$ of all complex numbers.

 By Lemma 2.2, the coefficient of $x_1^{n-1}\ldots x_n^{n-1}$ in the polynomial
 $$P(x_1,\ldots,x_n):=\prod_{1\ls i<j\ls n}(x_j-x_i)(jx_j-ix_i)\in \C[x_1,\ldots,x_n]$$
 is $(-1)^{n(n-1)/2}\per[i^{j-1}]_{1\ls i,j\ls n}$, which is nonzero since $\per[i^{j-1}]_{1\ls i,j\ls n}>0$.
  Applying Lemma 2.1 we see that there are $x_1,\ldots,x_n\in A=\{a_1,\ldots,a_n\}$
 with $P(x_1,\ldots,x_n)\not=0$. Thus, for some $\pi\in S_n$ all the numbers
 $ka_{\sigma(k)}\ (k=1,\ldots,n)$ are distinct. This ends the proof of part (i).

 (ii) Let $e$ be the identity of the group $G$.
 Suppose that $a=b^2$ and also $a^2=b$. Then $a=(a^2)^2=a^4$, and hence $a^3=e$. As the order of $a$ is not three, we have $a=e$ and hence $b=a^2=e$, which leads to a contradiction since $a\not=b$. Therefore $a^{\sigma(1)}$ and $b^{\sigma(2)}$ are distinct for some $\sigma\in S_2$.

To prove the second assertion in Theorem 1.4(ii), we distinguish two cases.
\smallskip

{\it Case} 1. One of $a,b,c$ is the square of another element among $a,b,c$.

Without loss of generality, we simply assume that $a=b^2$. As $a\not=b$ we have $b\not=e$.
As $b$ is not of order two, we also have $a\not=e$. Note that $b^2=a\not=c$.
If $b^2=a^3$, then $a=a^3$ which is impossible since the order of $a$ is not two.
If $a^3\not=c$, then $c,b^2,a^3$ are pairwise distinct.

Now assume that $a^3=c$. As $a$ is not of order three, we have $b\not=a^2$ and $c\not=e$.
Note that $a^3=c\not=b$ and also $a^3=c\not=c^2$. If $b\not=c^2$, then $b,c^2,a^3$ are pairwise distinct.
If $b=c^2$, then $a=b^2=c^4=(a^3)^4$ and hence the order of $a$ is $11$, thus $a^2\not=(a^3)^3=c^3$
and hence $b,a^2,c^3$ are pairwise distinct.
\smallskip

{\it Case} 2. None of $a,b,c$ is the square of another one among $a,b,c$.

Suppose that there is no $\tau\in S_3$ with $a^{\tau(1)},b^{\tau(2)},c^{\tau(3)}$ pairwise distinct.
Then $c^3\in\{a,b^2\}\cap\{a^2,b\}$. If $c^3=a$,
 then $c^3\not=b$ and hence $a=c^3=a^2$, thus $a=e=c$ which leads to a contradiction.
 (Recall that none of $a,b,c$ is of order $3$.)
 Therefore $c^3=b^2$. As $c$ is not of order three, if $b=e$ then we have $c=e=b$ which is impossible. So $c^3=b^2\not=b$
 and hence $b^2=c^3=a^2$. Similarly, $a^3=b^2=c^2$. Thus $a^3=b^2=a^2$, hence $a=e$ and $b^2=a^2=e$,
 which contradicts $b\not=a$ since $b$ is not of order two.

 In view of the above, we have finished the proof of Theorem 1.4. \qed

\medskip
\noindent{\it Proof of Theorem 1.5}. The subgroup of $G$ generated by $a_1,\ldots,a_n$ is a finitely generated torsion-free abelian group.
So we may simply assume that $G=\Z^r$ for some positive integer $r$ without any loss of generality.
It is well known that there is a linear ordering $\ls$ on $G=\Z^r$ such that for any $a,b,c\in G$
if $a<b$ then $-b<-a$ and $a+c<b+c$ (cf. [L]). For convenience, we suppose $a_1<a_2<\ldots<a_n$ without any loss of generality.

If $n=4$, then $(b_1,b_2,b_3,b_4)=(a_1,a_2,a_3,a_4)$ meets the requirement since
 $$a_1+a_2+a_3<a_4+a_1+a_2<a_3+a_4+a_1<a_2+a_3+a_4.$$
 Below we assume $n\gs5$.

Clearly
$$a_1+a_2+a_3<a_2+a_3+a_4<\ldots<a_{n-2}+a_{n-1}+a_n.$$
For convenience we set $$S:=\{a_{i-1}+a_i+a_{i+1}:\ i=2,\ldots,n-1\}, $$
and let $\min S$ and $\max S$ denote the least element and the largest element of $S$ respectively.
Note that
$$\min S=a_1+a_2+a_3<a_n+a_1+a_2<a_{n-1}+a_n+a_1<\max S=a_{n-2}+a_{n-1}+a_n.$$
If $\{a_n+a_1+a_2,\ a_{n-1}+a_n+a_1\}\cap S=\emptyset$, then
 $(b_1,\ldots,b_n)=(a_1,\ldots,a_n)$ meets the requirement.
Obviously
$$-a_n<-a_{n-1}<\ldots<-a_2<-a_1\ \ \t{and}\ \ (-a_2)+(-a_1)+(-a_n)=-(a_1+a_2+a_n).$$
So, it suffices to find a desired permutation $b_1,\ldots,b_n$ of $a_1,\ldots,a_n$ under the condition $a_{n-1}+a_n+a_1\in S$.

{\it Case} 1. $n=5$.

As $a_4+a_5+a_1\in S$, we have $a_4+a_5+a_1=a_2+a_3+a_4$ and
we may take $(b_1,\ldots,b_5)=(a_1,a_2,a_3,a_5,a_4)$ since
$$a_1+a_2+a_3<a_4+a_1+a_2<a_2+a_3+a_4=a_5+a_4+a_1<a_2+a_3+a_5<a_3+a_5+a_4.$$

{\it Case} 2. $n=6$.

As $a_5+a_6+a_1\in S$, the sum $a_5+a_6+a_1$ is equal to $a_2+a_3+a_4$ or $a_3+a_4+a_5$.
If $a_5+a_6+a_1=a_2+a_3+a_4$, then we may take
$(b_1,\ldots,b_6)=(a_1,a_2,a_5,a_3,a_4,a_6)$ since
$$\align a_1+a_2+a_5<&a_6+a_1+a_2<a_4+a_6+a_1<a_5+a_6+a_1=a_2+a_3+a_4
\\<&a_2+a_5+a_3<a_5+a_3+a_4<a_3+a_4+a_6.
\endalign$$
If $a_5+a_6+a_1=a_3+a_4+a_5$, then $a_6+a_1=a_3+a_4$ and we may take
$(b_1,\ldots,b_6)=(a_1,a_2,a_3,a_4,a_6,a_5)$ since
$$\align a_1+a_2+a_3<&a_5+a_1+a_2<a_6+a_1+a_2=a_2+a_3+a_4
\\<&a_3+a_4+a_5=a_6+a_5+a_1<a_3+a_4+a_6<a_4+a_6+a_5.
\endalign$$

{\it Case} 3. $n=7$.

As $a_6+a_7+a_1\in S$, the sum $a_6+a_7+a_1$ is equal to $a_2+a_3+a_4$ or $a_3+a_4+a_5$ or $a_4+a_5+a_6$.
If $a_6+a_7+a_1=a_4+a_5+a_6$, then $a_7+a_1=a_4+a_5$ and we may take
$(b_1,\ldots,b_7)=(a_2,a_1,a_4,a_5,a_3,a_6,a_7)$ since
$$\align a_2+a_1+a_4<&a_1+a_4+a_5=a_1+a_1+a_7<a_7+a_2+a_1
\\<&a_7+a_1+a_3=a_4+a_5+a_3<a_5+a_3+a_6
\\<&a_4+a_5+a_6=a_1+a_6+a_7<a_2+a_6+a_7<a_3+a_6+a_7.
\endalign$$
If $a_6+a_7+a_1=a_2+a_3+a_4$, then we may take
$(b_1,\ldots,b_7)=(a_1,a_2,a_3,a_5,a_4,a_6,a_7)$ since
$$\align a_1+a_2+a_3<&a_7+a_1+a_2<a_5+a_7+a_1<a_6+a_7+a_1=a_2+a_3+a_4
\\<&a_2+a_3+a_5<a_3+a_5+a_4<a_5+a_4+a_6<a_4+a_6+a_7.
\endalign$$
If $a_6+a_7+a_1=a_3+a_4+a_5$ and $a_5+a_6+a_1\not=a_2+a_3+a_4$, then $a_6+a_1<a_3+a_4$ and we may take
$(b_1,\ldots,b_7)=(a_1,a_2,a_3,a_4,a_7,a_5,a_6)$ since
$$\align a_1+a_2+a_3<&a_6+a_1+a_2<\min\{a_5+a_6+a_1,a_2+a_3+a_4\}
\\<&\max\{a_5+a_6+a_1,a_2+a_3+a_4\}<a_1+a_6+a_7=a_3+a_4+a_5
\\<&a_3+a_4+a_7<a_4+a_7+a_5<a_7+a_5+a_6.
\endalign$$
If $a_6+a_7+a_1=a_3+a_4+a_5$ and $a_5+a_6+a_1=a_2+a_3+a_4$, then $a_7+a_1<a_3+a_4$ and we may take
$(b_1,\ldots,b_7)=(a_1,a_2,a_3,a_4,a_6,a_5,a_7)$ since
$$\align a_1+a_2+a_3<&a_7+a_1+a_2<a_5+a_6+a_1=a_2+a_3+a_4
\\<&a_5+a_7+a_1<a_3+a_4+a_5=a_6+a_7+a_1
\\<&a_3+a_4+a_6<a_4+a_6+a_5<a_6+a_5+a_7.
\endalign$$

{\it Case} 4. $n>7$ and $a_n+a_1+a_2\not\in S$.

In this case, there is a unique $2<i<n-1$ with $a_{i-1}+a_i+a_{i+1}=a_{n-1}+a_n+a_1$.
If $i<n-3$, then we may take
$$(b_1,\ldots,b_n)=(a_1,\ldots,a_{i-2},a_{i-1},a_i,a_{i+2},a_{i+1},a_{i+3},\ldots,a_n)$$
because
$$\align a_{i-2}+a_{i-1}+a_i<&a_{i-1}+a_i+a_{i+1}=a_{n-1}+a_n+a_1<a_{i-1}+a_i+a_{i+2}
\\<&a_i+a_{i+2}+a_{i+1}<a_{i+2}+a_{i+1}+a_{i+3}
\\<&a_{i+1}+a_{i+3}+a_{i+4}<\ldots<a_{n-2}+a_{n-1}+a_n.
\endalign$$
When $i\in\{n-2,n-3\}$, we have $i\gs n-3>4$, and hence in the case $a_1+a_2+a_n\not=a_{i-4}+a_{i-3}+a_{i-1}$ we may take
$$(b_1,\ldots,b_n)=(a_1,\ldots,a_{i-4},a_{i-3},a_{i-1},a_{i-2},a_i,a_{i+1},a_{i+2},\ldots,a_n)$$
because
$$\align a_{i-4}+a_{i-3}+a_{i-2}<&a_{i-4}+a_{i-3}+a_{i-1}<a_{i-3}+a_{i-1}+a_{i-2}
\\<&a_{i-1}+a_{i-2}+a_i<a_{i-2}+a_i+a_{i+1}
\\<&a_{i-1}+a_i+a_{i+1}=a_{n-1}+a_n+a_1
\\<&a_i+a_{i+1}+a_{i+2}<\ldots<a_{n-2}+a_{n-1}+a_n
\endalign$$
and
$$\align a_n+a_1+a_2<&(a_{i-2}+a_{n-1}-a_{i+1})+a_n+a_1
\\&=a_{i-2}-a_{i+1}+(a_{i-1}+a_i+a_{i+1})=a_{i-1}+a_{i-2}+a_i.
\endalign$$
If $i\in\{n-2,n-3\}$ and $a_1+a_2+a_n=a_{i-4}+a_{i-3}+a_{i-1}$, then we may take
$$(b_1,\ldots,b_n)=(a_1,\ldots,a_{i-4},a_{i-3},a_i,a_{i-2},a_{i-1},a_{i+1},a_{i+2},\ldots,a_n)$$
because
$$\align a_n+a_1+a_2=&a_{i-4}+a_{i-3}+a_{i-1}
\\<&a_{i-4}+a_{i-3}+a_i<a_{i-3}+a_i+a_{i-2}<a_i+a_{i-2}+a_{i-1}
\\<&a_{i-2}+a_{i-1}+a_{i+1}<a_{i-1}+a_i+a_{i+1}=a_{n-1}+a_n+a_1
\\<&a_{i-1}+a_{i+1}+a_{i+2}<\ldots<a_{n-2}+a_{n-1}+a_n.
\endalign$$

{\it Case} 5. $n>7$ and $a_n+a_1+a_2\in S$.

In this case, for some $2<j<i\ls n-2$ we have
$$a_{n-1}+a_n+a_1=a_{i-1}+a_i+a_{i+1}>a_{j-1}+a_j+a_{j+1}=a_n+a_1+a_2.$$
If $j+1=i$, then
$$\align a_{n-1}-a_2=&(a_{n-1}+a_n+a_1)-(a_n+a_1+a_2)
\\=&a_{i-1}+a_i+a_{i+1}-(a_i+a_{i-1}+a_{i-2})=a_{i+1}-a_{i-2}
\endalign$$
which is impossible since $i\gs 4$ and $n>6$.

If $i-j>5$, then
$$(b_1,\ldots,b_n)=(a_1,\ldots,a_{j-1},a_j,a_{j+2},a_{j+1},a_{j+3},\ldots,a_{i-3},a_{i-1},a_{i-2},a_i,a_{i+1},\ldots,a_n)$$
meets the requirement since
$$\align a_{j-1}+a_j+a_{j+1}&=a_n+a_1+a_2<a_{j-1}+a_j+a_{j+2}
\\&<a_j+a_{j+2}+a_{j+1}<a_{j+2}+a_{j+1}+a_{j+3}
\\&<\ldots<a_{i-3}+a_{i-1}+a_{i-2}<a_{i-1}+a_{i-2}+a_i
\\&<a_{i-2}+a_i+a_{i+1}<a_{i-1}+a_i+a_{i+1}=a_{n-1}+a_n+a_1
\\&<a_i+a_{i+1}+a_{i+2}<\ldots<a_{n-2}+a_{n-1}+a_n.\endalign$$
If $i-j=5$, then $j+4=i-1$ and
$$(b_1,\ldots,b_n)=(a_1,\ldots,a_{j-1},a_j,a_{j+2},a_{j+1},a_{i-1},a_{i-2},a_i,a_{i+1},\ldots,a_n)$$
meets the requirement. If $i-j=4$, then
$$(b_1,\ldots,b_n)=(a_1,\ldots,a_{j-1},a_j,a_{j+2},a_{j+3},a_{j+1},a_i,a_{i+1},\ldots,a_n)$$
meets the requirement since
$$\align a_{j-1}+a_j+a_{j+1}=&a_n+a_1+a_2
\\<&a_{j-1}+a_j+a_{j+2}<a_j+a_{j+2}+a_{j+3}
\\<&a_{j+2}+a_{j+3}+a_{j+1}<a_{j+3}+a_{j+1}+a_i
\\<&a_{j+1}+a_i+a_{i+1}<a_{i-1}+a_i+a_{i+1}=a_{n-1}+a_n+a_1
\\<&a_i+a_{i+1}+a_{i+2}<\ldots<a_{n-2}+a_{n-1}+a_n.
\endalign$$
If $i-j=3$, then
$$(b_1,\ldots,b_n)=(a_1,\ldots,a_{j-1},a_j,a_{j+2},a_{j+1},a_i,a_{i+1},\ldots,a_n)$$
meets the requirement since
$$\align a_{j-1}+a_j+a_{j+1}=&a_n+a_1+a_2
\\<&a_{j-1}+a_j+a_{j+2}<a_j+a_{j+2}+a_{j+1}
\\<&a_{j+2}+a_{j+1}+a_i=a_{i-1}+a_{i-2}+a_i<a_{i-2}+a_i+a_{i+1}
\\<&a_{i-1}+a_i+a_{i+1}=a_{n-1}+a_n+a_1
\\<&a_i+a_{i+1}+a_{i+2}<\ldots<a_{n-2}+a_{n-1}+a_n.
\endalign$$
If $j>4$ and $i=j+2$, then
$$(b_1,\ldots,b_n)=(a_1,\ldots,a_{j-3},a_{j-1},a_{j-2},a_{j+1},a_j,a_i,a_{i+1},a_{i+2},\ldots,a_n)$$
meets the requirement since
$$\align a_{j-4}+a_{j-3}+a_{j-1}&<a_{j-3}+a_{j-1}+a_{j-2}<a_{j-1}+a_{j-2}+a_{j+1}
\\&<a_{j-2}+a_{j+1}+a_j<a_{j-1}+a_j+a_{j+1}=a_n+a_1+a_2
\\&<a_{j+1}+a_j+a_i<a_j+a_i+a_{i+1}
\\&<a_{i-1}+a_i+a_{i+1}=a_{n-1}+a_n+a_1<a_i+a_{i+1}+a_{i+2}.\endalign$$
If $i=j+2\ls n-4$, then
$$(b_1,\ldots,b_n)=(a_1,\ldots,a_{j-2},a_{j-1},a_j,a_i,a_{i-1},a_{i+2},a_{i+1},a_{i+3},a_{i+4},\ldots,a_n)$$
meets the requirement since
$$\align a_{j-2}+a_{j-1}+a_j&<a_{j-1}+a_j+a_{j+1}=a_n+a_1+a_2
\\&<a_{j-1}+a_j+a_i<a_j+a_i+a_{i-1}
\\&<a_{i-1}+a_i+a_{i+1}=a_{n-1}+a_n+a_1
\\&<a_i+a_{i-1}+a_{i+2}<a_{i-1}+a_{i+2}+a_{i+1}
\\&<a_{i+2}+a_{i+1}+a_{i+3}<a_{i+1}+a_{i+3}+a_{i+4}
\\&<\ldots<a_{n-2}+a_{n-1}+a_n.\endalign$$
If $i\gs n-3$, $j\ls 4$ and $i-j=2$, then $2=i-j\gs n-3-4$ and hence $n\in\{8,9\}$.

For $n=8$, we need to consider the case $i=6$ and $j=4$.
As
$a_8+a_1+a_2=a_3+a_4+a_5$ and $a_7+a_8+a_1=a_5+a_6+a_7$,
we have $a_8+a_1=a_5+a_6=a_3+a_4+a_5-a_2$.
If $2a_5\not=a_4+a_7$, then $a_5+a_8+a_1=2a_5+a_6\not=a_4+a_6+a_7$ and hence we may take
$$(b_1,\ldots,b_8)=(a_1,a_2,a_3,a_4,a_6,a_7,a_5,a_8)$$
since
$$\align a_1+a_2+a_3&<a_2+a_3+a_4<a_3+a_4+a_5=a_8+a_1+a_2<a_3+a_4+a_6
\\&<\min\{a_4+a_6+a_7,a_5+a_8+a_1\}<\max\{a_4+a_6+a_7,a_5+a_8+a_1\}
\\&<a_6+a_7+a_5=a_7+a_8+a_1<a_7+a_5+a_8.
\endalign$$
If $2a_5=a_4+a_7$, then $a_6+a_8+a_1=a_5+2a_6>a_4+a_5+a_7$ and we may take
$$(b_1,\ldots,b_8)=(a_1,a_2,a_3,a_4,a_5,a_7,a_8,a_6)$$
since
$$\align a_1+a_2+a_3&<a_1+a_3+a_4=a_1+a_2+a_6<a_2+a_3+a_4
\\&<a_3+a_4+a_5=a_8+a_1+a_2<a_4+a_5+a_7<a_6+a_8+a_1
\\&<a_5+a_7+a_8<a_7+a_8+a_6.
\endalign$$
When $n=8$, $i=5$ and $j=3$, it suffices to apply the result for $i=6$ and $j=4$ to the sequence
$$\align a_1'=-a_8<&a_2'=-a_7<a_3'=-a_6<a_4'=-a_5
\\<&a_5'=-a_4<a_6'=-a_3<a_7'=-a_2<a_8'=-a_1
\endalign$$
since $a_7'+a_8'+a_1'=-(a_1+a_2+a_8)=-(a_2+a_3+a_4)=a_5'+a_6'+a_7'$ and
$a_8'+a_1'+a_2'=-(a_1+a_7+a_8)=-(a_4+a_5+a_6)=a_3'+a_4'+a_5'$.

Now it remains to consider the last case where $n=9$, $i=6$ and $j=4$.
As $a_3+a_4+a_5=a_9+a_1+a_2$ and $a_5+a_6+a_7=a_8+a_9+a_1$, we have
$a_3+a_4<a_9+a_1$ and hence $a_3+a_4+a_6<a_3+a_4+a_7<a_7+a_9+a_1$.
If $a_7+a_9+a_1=a_4+a_5+a_6$, then
$$a_8-a_7=(a_8+a_9+a_1)-(a_7+a_9+a_1)=a_5+a_6+a_7-(a_4+a_5+a_6)=a_7-a_4.$$
When $2a_7\not=a_8+a_4$, we have $a_7+a_9+a_1\not=a_4+a_5+a_6$ and hence we may take
$$(b_1,\ldots,b_9)=(a_1,a_2,a_3,a_4,a_6,a_5,a_8,a_7,a_9)$$
since
$$\align a_1+a_2+a_3&<a_2+a_3+a_4<a_3+a_4+a_5=a_9+a_1+a_2<a_3+a_4+a_6
\\&<\min\{a_4+a_5+a_6,a_7+a_9+a_1\}<\max\{a_4+a_5+a_6,a_7+a_9+a_1\}
\\&<a_6+a_5+a_7=a_8+a_9+a_1<a_6+a_5+a_8
\\&<a_5+a_8+a_7<a_8+a_7+a_9.
\endalign$$
If $2a_7=a_8+a_4$, then $a_5+a_6+a_7<2a_7+a_6=a_4+a_6+a_8$ and hence
we may take
$$(b_1,\ldots,b_9)=(a_1,a_2,a_3,a_4,a_6,a_8,a_5,a_7,a_9)$$
since
$$\align a_1+a_2+a_3&<a_2+a_3+a_4<a_3+a_4+a_5=a_9+a_1+a_2<a_3+a_4+a_6
\\&<a_9+a_1+a_6<a_7+a_9+a_1<a_8+a_9+a_1=a_5+a_6+a_7
\\&<a_4+a_6+a_8<a_6+a_8+a_5<a_8+a_5+a_7<a_5+a_7+a_9.
\endalign$$

In view of the above, we have completed the proof of Theorem 1.5. \qed

\heading{4. Some conjectures}\endheading

Motivated by Theorem 1.3(i) and Theorem 1.4, we pose the following conjecture for finite groups.

\proclaim{Conjecture 4.1} Let $n$ be a positive integer, and let $G$ be a group
containing no element of order among $2,\ldots,n+1$.
Then, for any $A\se G$ with $|A|=n$, we may write $A=\{a_1,\ldots,a_n\}$
with $a_1,a_2^2,\ldots,a_n^n$ pairwise distinct.
\endproclaim
\Remark\ 4.1. (a) Theorem 1.4 shows that this conjecture holds when $n\ls3$ or $G$ is a torsion-free abelian group.

(b) For $n=4,5,6,7,8,9$ we have verified the conjecture for cyclic groups $G=\Z/m\Z$ with $|G|=m$ not exceeding
$100,\, 100,\, 70,\, 60,\, 30,\, 30$ respectively.

(c) If $G$ is a finite group with $|G|>1$, then
the least order of a non-identity element of $G$ is $p(G)$, the smallest prime divisor of $|G|$.
\medskip

Inspired by Theorem 1.3, we formulate the following conjecture.

\proclaim{Conjecture 4.2} Let $n>1$ be an integer with $n\not\eq2\pmod4$.

{\rm (i)} We have
$$\per[i^{j-1}]_{1\ls i,j\ls n-1}\eq0\pmod n.\tag4.1$$

{\rm (ii)} If $n\eq1\pmod3$, then
$$\per[i^{j-1}]_{1\ls i,j\ls n-1}\eq0\pmod {n^2}.\tag4.2$$
\endproclaim
\Remark\ 4.2. We have checked this conjecture via computing $\per[i^{j-1}]_{1\ls i,j\ls n-1}$ modulo $n^2$ for $n\ls 17$.
The sequence $a_n=\per[i^{j-1}]_{1\ls i,j\ls n}\ (n=1,2,3,\ldots)$ is available from [S18, A322363].
\medskip

\proclaim{Conjecture 4.3} {\rm (i)} For any $n\in\Z^+$, there is a permutation $\sigma_n\in S_n$
such that $k\sigma_n(k)+1$ is prime for every $k=1,\ldots,n$.

{\rm (ii)} For any integer $n>2$, there is a permutation $\tau_n\in S_n$ such that $k\tau_n(k)-1$ is prime for every $k=1,\ldots,n$.
\endproclaim
\Remark\ 4.3. See [S18, A321597] for related data and examples.

\proclaim{Conjecture 4.4} {\rm (i)} For each $n\in\Z^+$, there is a permutation $\pi_n$ of $\{1,\ldots,n\}$ such that $k^2+k\pi_n(k)+\pi_n(k)^2$ is prime for every $k=1,\ldots,n$.

{\rm (ii)} For any positive integer $n\not=7$, there is a permutation $\pi_n$ of $\{1,\ldots,n\}$ such that $k^2+\pi_n(k)^2$ is prime for every $k=1,\ldots,n$.
\endproclaim
\Remark\ 4.4. See [S18, A321610] for related data and examples.
\medskip

As usual, for $k=1,2,3,\ldots$ we let $p_k$ denote the $k$-th prime.

\proclaim{Conjecture 4.5} For any $n\in\Z^+$, there is a permutation $\pi\in S_n$
such that $p_k+p_{\pi(k)}+1$ is prime for every $k=1,\ldots,n$.
\endproclaim
\Remark\ 4.5. See [S18, A321727] for related data and examples.
\medskip

In 1973 J.-R. Chen [Ch] proved that there are infinitely many primes $p$
with $p+2$ a product of at most two primes; nowadays such primes $p$ are called Chen primes.

\proclaim{Conjecture 4.6} Let $n\in\Z^+$. Then,
there is an even permutation $\sigma\in S_n$ with $p_kp_{\sigma(k)}-2$ prime for all $k = 1,\ldots,n$. If $n>2$, then there is an odd permutation $\tau\in S_n$ with $p_kp_{\tau(k)}-2$ prime for all $k =1,\ldots,n$.
\endproclaim
\Remark\ 4.6. See [S18, A321855] for related data and examples. If we let $b(n)$ denote the number of even permutations $\sigma\in S_n$ with $p_kp_{\sigma(k)}-2$ prime for all $k = 1,\ldots,n$,  then
$$(b(1),\ldots,b(11)) = (1,1,1,1,3,6,1,1,33,125,226).$$
Conjecture 2.17(ii) of Sun [S15] implies that for any odd integer $n>1$ there is a prime $p\ls n$ such that $pn-2$ is prime.
\medskip

In 2002, Cloitre [Cl, A073112] created the sequence A073112 on OEIS whose $n$-th term
is the number of permutations $\pi\in S_n$ with $\sum_{k=1}^n\f1{k+\pi(k)}\in\Z$.
Recently Sun [S18a] conjectured that for any integer $n>5$
there is a permutation $\pi\in S_n$ satisfying
$$\sum_{k=1}^n\f1{k+\pi(k)}=1,$$ and this was later confirmed
by the user Zhao Shen at Mathoverflow via clever induction arguments.

In 1982 A. Filz (cf. [G, pp. 160-162]) conjectured that for any $n=2,4,6,\ldots$ there is a circular permutation
$(i_1,\ldots,i_n)$ of $1,\ldots,n$ such that all the $n$ adjacent sums
$$i_1+i_2,\ i_2+i_3,\ \ldots,\ i_{n-1}+i_n,\ i_n+i_1$$ are prime.

Motivated by this, we pose the following conjecture.

\proclaim{Conjecture 4.7}
{\rm (i)} For any integer $n>6$, there is a permutation $\pi\in S_n$ such that
$$\sum_{k=1}^{n-1}\f1{\pi(k)+\pi(k+1)}=1.\tag4.3$$
Also, for any integer $n>7$, there is a permutation $\pi\in S_n$ such that
$$\f1{\pi(1)+\pi(2)}+\f1{\pi(2)+\pi(3)}+\ldots+\f1{\pi(n-1)+\pi(n)}+\f1{\pi(n)+\pi(1)}=1.\tag4.4$$

{\rm (ii)} For any integer $n>7$, there is a permutation $\pi\in S_n$ such that
$$\sum_{k=1}^{n-1}\f1{\pi(k)^2-\pi(k+1)^2}=0.\tag4.5$$
\endproclaim
\Remark\ 4.7. See [S18, A322070 and A322099] for related data and examples.
For the latter assertion in Conjecture 4.7(i), the equality (4.4) with $n=8$ holds
if we take $(\pi(1),\ldots,\pi(8))=(6,1,5,2,4,3,7,8)$.
In a previous version of this paper posted to arXiv, the author also conjectured that
for any integer $n>5$ there is a permutation $\pi\in S_n$ with
$\sum_{k=1}^{n-1}\f1{\pi(k)\pi(k+1)}=1$; this, together with two other conjectures of the author,
was confirmed by G.-N. Han \cite{H}.

\proclaim{Conjecture 4.8} {\rm (i)} For any integer $n>1$, there is a permutation $\pi\in S_n$ such that
$$\sum_{0<k<n}\pi(k)\pi(k+1)\in\{2^m+1:\ m=0,1,2,\ldots\}.\tag4.6$$

{\rm (ii)} For any integer $n>4$, there is a unique power of two which can be written as $\sum_{k=1}^{n-1}\pi(k)\pi(k+1)$ with $\pi\in S_n$ and $\pi(n)=n$.
\endproclaim
\Remark\ 4.8. Concerning part (i) of Conjecture 4.8, when $n=4$ we may choose $(\pi(1),\ldots,\pi(4))=(1,3,2,4)$ so that
$$\sum_{k=1}^3\pi(k)\pi(k+1)=1\times3+3\times2+2\times4=2^4+1.$$

For any $\pi\in S_n$, if for each $k=1,\ldots,n$ we let
$$\pi'(k)=\cases\pi(\pi^{-1}(k)+1)&\t{if}\ \pi^{-1}(k)\not=n,
\\\pi(1)&\t{if}\ \pi^{-1}(k)=n,\endcases$$
then $\pi'\in S_n$ and
$$\pi(1)\pi(2)+\ldots+\pi(n-1)\pi(n)+\pi(n)\pi(1)=\sum_{k=1}^nk\pi'(k).$$

By the Cauchy-Schwarz inequality (cf. [N, p.\,178]), for any $\pi\in S_n$ we have
$$\bigg(\sum_{k=1}^nk\pi(k)\bigg)^2\ls\bigg(\sum_{k=1}^nk^2\bigg)\bigg(\sum_{k=1}^n\pi(k)^2\bigg)$$
and hence
$$\sum_{k=1}^nk\pi(k)\ls\sum_{k=1}^nk^2=\frac{n(n+1)(2n+1)}6.$$
If we let $\sigma(k)=n+1-\pi(k)$ for all $k=1,\ldots,n$, then $\sigma\in S_n$ and
$$\align\sum_{k=1}^n k\pi(k)=&\sum_{k=1}^nk(n+1-\sigma(k))=(n+1)\sum_{k=1}^nk-\sum_{k=1}^nk\sigma(k)
\\\gs&\frac{n(n+1)^2}2-\frac{n(n+1)(2n+1)}6=\frac{n(n+1)(n+2)}6.
\endalign$$
Thus
$$\bg\{\sum_{k=1}^nk\pi(k):\ \pi\in S_n\bg\}\subseteq T(n):=\left\{\frac{n(n+1)(n+2)}6,\ldots,\frac{n(n+1)(2n+1)}6\right\}.\tag4.7$$
Actually equality in (4.7) holds when $n\not=3$, which was first realized by M. Aleksevev (cf. the comments in [B]). Note that
$|T(n)|=n(n^2-1)/6+1$.

Inspired by the above analysis, here we pose the following conjecture in additive combinatorics.

\proclaim{Conjecture 4.9} Let $n\in\Z^+$ and let $F$ be a field with $p(F)>n+1$, where $p(F)=p$ if the characteristic of $F$ is a prime $p$, and $p(F)=+\infty$ if the characteristic of $F$ is zero. Let $A$ be any finite subset of $F$ with $|A|\gs n+\delta_{n,3}$, where $\delta_{n,3}$ is $1$ or $0$ according as $n=3$ or not. Then, for the set
$$S(A):=\bigg\{\sum_{k=1}^n ka_k:\ a_1,\ldots,a_n\ \text{are distinct elements of}\ A\bigg\},\tag4.8$$
we have $$|S(A)|\gs\min\left\{p(F),\ (|A|-n)\frac{n(n+1)}2+\frac{n(n^2-1)}6+1\right\}.\tag4.9$$
\endproclaim
\Remark\ 4.9. One may compare this conjecture with the author's
conjectural linear extension of the Erd\H os-Heilbronn conjecture (cf. [SZ]).
Perhaps, Conjecture 4.9 remains valid if we replace the field $F$ by
any finite additive group $G$ with $|G|>1$
and use $p(G)$ (the least prime factor of $|G|$) instead of $p(F)$.
\medskip

Recall that the torsion subgroup of a group $G$ is given by
$$\t{Tor}(G)=\{g\in G:\ g\ \t{is of finite order}\}.$$

Conjecture 3.3(i) of the author \cite{S19}
states that if $A$ is an $n$-subset (with $|A|=n>2$) of an additive abelian group $G$ of odd order
then there is a numbering $a_1,\ldots,a_n$ of all the elements of $A$ such that $a_1+a_2,\ldots,a_{n-1}+a_n,a_n+a_1$ are pairwise distinct, this was verified by Yu-Xuan Ji
(a student at Nanjing Univ.) for $|G|<30$ in 2020.
Motivated by this and Theorem 1.5, we formulate the following conjecture.

\proclaim{Conjecture 4.10} Let $G$ be an additive abelian group with $\Tor(G)$
cyclic or $|\Tor(G)|$ odd. For any finite subset $A$ of $G$ with $|A|=n>3$,
there is a numbering $a_1,\ldots,a_n$
of all the elements of $A$ such that the $n$ sums
$$a_1+a_2+a_3,\ a_2+a_3+a_4,\ \ldots,\ a_{n-2}+a_{n-1}+a_n,\ a_{n-1}+a_n+a_1,\ a_n+a_1+a_2$$
are pairwise distinct.
\endproclaim
\Remark\ 4.10. (a) Conjecture 4.10 holds in the case $A=G=\Z/n\Z=\{\bar a=a+n\Z:\ a\in\Z\}$ with $n>3$ and $3\nmid n$ since the natural list $\bar 0,\bar 1,\ldots,\overline{n-1}$
of the elements of $\Z/n\Z$ meets the requirement.

 (b) In 2008 the author [S08] proved that for any three $n$-subsets $A,B,C$ of an additive abelian group $G$ with $\Tor(G)$ cyclic, there
is a numbering $a_1,\ldots,a_n$ of the elements of $A$, a numbering $b_1,\ldots,b_n$ of the elements of $B$ and a numbering $c_1,\ldots,c_n$ of the elements of $C$ such that the $n$ sums $a_1+b_1+c_1,\ldots,a_n+b_n+c_n$ are pairwise distinct.


\medskip

\widestnumber\key{DKSS}

 \Refs

\ref\key A\by N. Alon\paper Combinatorial Nullstellensatz\jour
Combin. Probab. Comput. \vol \yr 1999\pages 7--29\endref

\ref\key B\by J. Boscole\paper {\rm Sequence A126972 in OEIS, 2007}
\jour Website: {\tt http://oeis.org/A126972}\endref

\ref\key Br\by P. Bradley\paper Prime number sums\jour preprint, arXiv:1809.01012 (2018)\endref

\ref\key Ch\by J.-R. Chen\paper On the representation of a larger even integer as the sum of a prime and the product of at most two primes\jour Sci. Sinica 16 (1973), 157--176\endref

\ref\key Cl\by B. Cloitre\paper {\rm Sequences A073112 and A073364 in OEIS (2002)}\jour {\tt http://oeis.org}\endref

\ref\key DKSS\by  S. Dasgupta, G. K\'arolyi, O. Serra and B. Szegedy\paper Transversals of additive
Latin squares\jour Israel J. Math. \vol 126 \yr 2001\pages 17¨C28\endref

\ref\key FSX\by T. Feng, Z.-W. Sun and Q. Xiang\paper Exterior algebras and two conjectures on finite abelian groups\jour Israel J. Math.\vol 182\yr 2011\pages 425--437\endref

\ref\key GS\by F. Ge and Z.-W. Sun\paper On a permutation problem for finite abelian groups
\jour Electron. J. Combin. \vol 24\yr 2017\pages no.\,1, \#P1.17, 1--6\endref

\ref\key G\by R. K. Guy\book {\rm Unsolved Problems in Number Theory}\publ 3rd Edition, Springer, 2004\endref

\ref\key H\by G.-N. Han\paper On the existence of permutations conditioned by certain rational functions
\jour Electron. Res. Arch. \vol 28\yr 2020\pages 149--156\endref

\ref\key He\by E. Hecke\book Lectures on the Theory of Algebraic Numbers, Grad. Texts in
Math., 77, Springer, New York, 1981, pp. 108¨C116.\endref

\ref\key L\by  F. W. Levi\paper Ordered groups\jour Proc. Indian Acad. Sci. Sect. A
\vol 16 \yr 1942\pages 256--263\endref

\ref\key N\by M. B. Nathanson\book {\rm Additive Number Theory: The Classical Bases, Grad. Texts in
Math., Vol. 164, Springer, New York, 1996}\endref

\ref\key S08\by Z.-W. Sun\paper An additive theorem and restricted sumsets
\jour Math. Res. Lett.\vol 15\yr 2008\pages 1263--1276\endref

\ref\key S15\by Z.-W. Sun \paper Problems on combinatorial properties of primes\jour
¡¡¡¡in: M. Kaneko, S. Kanemitsu and J. Liu (eds.), Number Theory: Plowing and Starring through High Wave Forms, Proc. 7th China-Japan Seminar
¡¡¡¡(Fukuoka, Oct. 28--Nov. 1, 2013), Ser. Number Theory Appl., Vol. 11, World Sci., Singapore, 2015, pp. 169--187\endref

\ref\key S18\by Z.-W. Sun\paper {\rm Sequences A321597, A321610, A321611, A321727,
A321855, A322070, A322099, A322363 in OEIS (2018)}
\jour {\tt http://oeis.org}\endref

\ref\key S18a\by Z.-W. Sun\paper Permutations $\pi\in S_n$ with $\sum_{k=1}^n\f1{k+\pi(k)}=1$
\jour Question 315648 on Mathoverflow, Nov. 19, 2018. Website: {\tt https://mathoverflow.net/questions/315648}\endref

\ref\key S19\by Z.-W. Sun\paper Some new problems in additive combinatorics
\jour Nanjing Univ. J. Math. Biquarterly \vol 36\yr 2019\pages 134--155.
{\tt http://maths.nju.edu.cn/$\sim$zwsun/196a.pdf}\endref

\ref\key SZ\by Z.-W. Sun and L.-L. Zhao\paper Linear extension of the Erd\H os-Heilbronn conjecture
\jour J. Combin. Theory Ser. A\vol 119\yr 2012\pages364--381\endref

\endRefs
\enddocument